\documentclass[12pt,reqno,sumlimits]{amsart}
\usepackage{amssymb}

\theoremstyle{plain}
\newtheorem{thm}{Theorem}[section]
\newtheorem{cor}[thm]{Corollary}
\newtheorem{lemma}[thm]{Lemma}
\newtheorem{prop}[thm]{Proposition}
\theoremstyle{definition}
\newtheorem{definition}[thm]{Definition}
\newtheorem{rmk}[thm]{Remark}

\newcommand{\R}{{\mathbb R}}

\newcommand{\C}{{\mathbb C}}
\newcommand{\V}{{\mathbb V}}

\newcommand{\g}{{\mathfrak  g}}

\newcommand{\calO}{{\mathcal O}}

\newcommand{\x}{{\mathfrak X}}
\newcommand{\y}{{\mathfrak Y}}
\renewcommand{\to}{\longrightarrow}
\newcommand{\dd}[1]{\frac{\partial}{\partial #1}}
\newcommand{\DD}[1]{\frac{d}{d #1}}

\newcommand{\Hom}{\operatorname{Hom}}

\newcommand{\Ad}{\operatorname{Ad}}

\newcommand{\Diff}{\operatorname{Diff}}

\newcommand{\spann}{{\operatorname{span}}}

\newcommand{\inv}{^{-1}}

  \numberwithin{equation}{section}

\begin{document}

\title{Virasoro Actions and Harmonic Maps (after Schwarz)}
\author{Mihaela Vajiac and Karen Uhlenbeck}
\address{Department of Mathematics\\
  University of Texas at Austin}
\email{mbvajiac@chapman.edu,uhlen@math.utexas.edu}

\maketitle

\section{Introduction}
\label{Introd}

The influence on theoretical physics on geometry and topology at this point
in time is overwhelming.  A great many mathematicians are working to verify
conjectures made by physicists or suggested by physics.  However, the most
difficult problem facing mathematicians is to clarify or translate the
intuition from quantum field theory or string theory which lies behind the
ideas.  In this article, we interpret a series of papers of by John Schwarz, a
leading originator and proponent of string theory, on Virasoro actions from
a decade ago [4,5]  We hope to shed some light on the mathematical origins of
Virasoro actions constructed in physics.

Our main result is that there is (formally) an infinitesimal action of a
complex half-Virasoro algebra on the space of harmonic maps from a simply
connected domain in $\C\cup(\inf)$ to the Lie group $SU(n)$.  This action is
an example of a family of actions defined on integrable systems.  The
action on KdV is probably the best known of these[8].  We spend some time
describing the origin of these type of actions. They occur in the context
of a loop group which is split via Riemann-Hilbert factorization. The
half-Virasoro algebra  acts infinitesimally on the group and restricts to
one factor   of the splitting. The action we are interested in is the
derived action on the   second factor.  We give a description in terms of
groups which is quite  transparent and does not involve formulae. This
description works well in the context of harmonic maps, and leads to the
more complicated formulas in the Lie algebra  setting.

The plan of the paper is as follows.  After this introduction, we review
the structural background on harmonic maps into a Lie group which originates
in a  paper of one of us [7] and can be found in the text by Guest [2].
In chapter 3 we review the definitions of triples of Lie groups and algebras
and describe the loop groups factorization used for harmonic maps.
Section 4 shows how the half-Virasoro actions arise in these factorizations.
The main results in the paper are in 5.1 and the  following corollaries.
Section 6   outlines the results in the Wick rotated version of harmonic
maps from $\R^{1,1}$ into $SU(n)$ and makes contact with the formulae of
Schwarz.

This is a small part of a project on Virasoro actions on integrable systems
which is joint with Chuu-Lian Terng.  We thank Dan Freed for much needed
inspiration and encouragement. We apologize for our non-inclusive
reference list.  The literature on harmonic maps, integrable systems and
Virasoro actions is immense and quite splintered.

Our results are incomplete, since to obtain the full Virasoro action
and the coupling to gravity proposed by Schwarz [4], it is necessary to extend
the actions to include the $L_{-1}$ generators.  From comparison to the
Virasoro actions on other integrable systems, these generators include
second flows in their description, if we regard the harmonic maps as
first flows.  We hope to extend the description in this direction, as well
as to treat the  extension to harmonic maps in other contexts.

\section{Background}
\label{sec:1}

We give a brief description following Guest~\cite{G1} to Lie-group
valued harmonic maps.

\begin{definition}
  \label{2.1}
  A {\em harmonic map} $s:\Omega\to G$, (where $\Omega$ is a simply connected
  domain in $\C\cup\{\infty\}$ and $G$ is a matrix Lie group,
  is a solution to the Euler-Lagrange equation:
  $$
  \frac{\partial}{\partial x}\left( s^{-1}\frac{\partial}{\partial x}\right)+
  \frac{\partial}{\partial y}\left(s^{-1}\frac{\partial}{\partial y}\right)=0.
  $$
  
  The map $s$ satisfies the reality condition $s^{-1}(q)=s^*(q),$ if
  $G=SU(n)$ .
\end{definition}
Let $L(SU(n))=\{s:\Omega\to SU(n),\, s \,{\text{ harmonic}} \}$.  The
Euler-Lagrange equations are equivalent to the following:
$$
(s^{-1}s_{\bar z})_z +(s^{-1}s_z)_{\bar z}=0
$$
Write $A=s^{-1}s_z$ and $B=s^{-1}s_{\bar z}$, then the harmonic map
equation becomes $A_{\bar z}+B_z=0,$ where $A,B:\C\to g\otimes\C$, and
$B=-A^*$.  An equivalent description becomes:
\begin{prop}
  \label{2.2}
  The harmonic map equation is equivalent to the system:
  \begin{eqnarray*}
    A_{\bar z}+B_z&=&0\\
    A_{\bar z}-B_z&=&[A,B],
  \end{eqnarray*}
  where $A,B:\C\to g\otimes\C,$ and $B+A^*=0.$
\end{prop}

The generalized solution associated to the harmonic map equation can
be constructed as follows.  Let $\lambda\in\C$,
$A_{\lambda}=\frac{1}{2}(1-\lambda^{-1})A$, $B_{\lambda}=\frac{1}{2}(1-\lambda)B,$ and
$B=-A^*.$ With this notation the harmonic map equation is equivalent
to the equation:
$$
(A_{\lambda})_{\bar z}-(B_{\lambda})_z=[A_{\lambda},B_{\lambda}], \qquad \forall \lambda\in\C-\{0\}.
$$

Note that $\lambda$ is a spectral (or twistor) parameter and should be
carefully distinguished from the spatial parameters $z,\, \bar{z}$.
It is easy to see that this represents the flatness condition of the
associated connections $D_\lambda=\{d+A(\lambda)\}$, which admit a frame of flat
sections. Let $E_\lambda$ be this frame.

Then, the harmonic map equations are equivalent to the following
system for the flat frame $E_\lambda$ of the associated connection:


\begin{eqnarray}
  \label{1}
  {E_\lambda}^{-1}\bar\partial E_\lambda &=&{E_\lambda}^{-1}(E_{\lambda })_z = \frac{1}{2}(1-\lambda^{-1})A,\nonumber\\
  {E_\lambda}^{-1}\partial E_\lambda &=&{E_\lambda}^{-1}(E_{\lambda })_{\bar z} = \frac{1}{2}(1-\lambda)B.
\end{eqnarray}

The harmonic map can be easily reconstructed from the flat frame
$E_{\lambda}$ as $s(z)=E_{-1}(z)$. More precisely, for the case of $SU(n)$,
which we treat in the rest of the paper, we have:

\begin{thm}[\cite{Uhlen1}]
  \label{2.3}
  If $s$ is harmonic and $s(p)\equiv I$, then there exists a unique
  $E:\C^*\times\Omega\to SL(n,\C)$ satisfying equations~\eqref{1} with
  \begin{itemize}
    \item[(a)] $E_1\equiv I,$
    \item[(b)] $E_{-1}=s,$
    \item[(c)] $E_\lambda(p)=I$.
    \item[(d)] $E_{\lambda^{-1}}^{-1}=E^*_{\bar{\lambda}}$.
  \end{itemize}
  Moreover, $E$ is analytic and holomorphic in $\lambda\in\C^*$. Note that
  $E_\lambda$ is unitary for $|\lambda|=1.$
\end{thm}

\begin{thm}[\cite{Uhlen1}] 
  \label{2.4}
  Suppose $E:\C^*\times\Omega\to SL(n,\C)$ is analytic
  and holomorphic in the first variable, satisfying the reality
  condition $E_{\lambda^{-1}}^{-1}=E^*_{\bar{\lambda}},\,$ $E_1\equiv I, E(p)=I$ and
  the expressions
  $$
  \frac{E_\lambda^{-1}(E_\lambda)_{\bar z}}{1-\lambda},\qquad
  \frac{E_\lambda^{-1}(E_\lambda)_{z}}{1-\lambda^{-1}}
  $$
  are constant in $\lambda$. Then $s=E_{-1}$ is harmonic.
\end{thm}

Thus, from the harmonic map $s$ we have obtained the ``extended
solution'' $E:\Omega\to\Omega G$, where $\Omega G$ is the loop group of $G$.

\begin{rmk}
  The choice of basepoint affects the extended solution, and hence the
  Virasoro actions. It is standard in integrable systems for the
  coice of basepoint to affect the constructions.

\end{rmk}


\section{Manin triples of groups and Riemann-Hilbert Factorization}


Let $(X^+,X^-,X)$ be a triple of Lie groups with $X^{\pm}\subset X$ and
$\mu:X^+\times X^-\to X$ a diffeomorphism, where $\mu(s_+,s_-)=s_+s_-.$ Then
$(X^+,X^-,X)$ is a {\em (Manin) triple of Lie groups}.  We say that
$(X^+,X^-,X)$ is a {\em local (Manin) triple} if 
$$
\mu:X^+\times X^-\simeq \tilde{X}\subset X
$$
is a diffeomorphism onto an open dense subset $\tilde{X}\subset X.$

We define projections 
\begin{eqnarray}
  P^{\pm}(s)&=&(\mu^{-1}(s))_{\pm},\nonumber\\
  P^{\pm}&:&\tilde{X}\to X^{\pm},
  \label{projections}
\end{eqnarray}
from the open set $\tilde{X}$ onto the subgroups $X^{\pm}.$

The standard example of a Manin triple is
$$(X^+,X^-,X)=(SU(n),\Delta(n),SL(n,\C)),$$ where $\Delta(n)$ is the group of upper
triangular matrices with real diagonal elements. The projection
$P^+:SL(n,\C)\to SU(n)$ is realized by applying the Gramm-Schmidt
process to the columns of a matrix.

In general, the projections $P^{\pm}$ do not have nice formulae.
However, at the Lie algebra level, the projection operators are easily
described. The triple of groups $(X^+,X^-,X)$ gives rise to a triple
of Lie algebras $(\x^+,\x^-,\x).$ At the Lie algebra level
$\x^{\pm}\subset\x$ and $\x=\x^++\x^-.$ The projection operators
$\Pi^{\pm}:\x\to\x^{\pm}$ exist everywhere, even when the group triple is
only local.  Note that
\begin{equation}
  \label{2}
  dP_{s^+}(s^+V)=s^+\Pi^+V,
\end{equation}  
for $s^+\in X^+, V\in\x.$ These are infinitesimal formulae which appear
later in the paper.


The examples of interest to us are local triples of groups where the
projections $P^{\pm}$ are realized by the Riemann-Hilbert
factorization.  We refer to Guest~\cite{G1} and
Pressley-Segal~\cite{PS} for a more detailed analysis of the sketch we
give of such factorizations.

For the general setting, we have a contour (not necessarily connected)
$\Gamma\subset S^2=\C\cup\{\infty\}$ and a sequence of open sets
$\calO_\epsilon^{\pm}\subset\calO_\delta^{\pm}$ for $0<\delta<\epsilon$, where

\begin{eqnarray*}
  S^2&\subset&\calO_\epsilon^{+}\cup\calO_\delta^{-},\quad \forall\epsilon,\delta\\
  \Gamma&=&\bigcap_\epsilon \calO_\epsilon^{+}\cap\calO_\epsilon^{-}.
\end{eqnarray*}
Regard $\calO_\epsilon^{-}$ as a thickening of the interior of $\Gamma$ and
$\calO_\epsilon^{+}$ as a thickening of the exterior. We define:
\begin{eqnarray*}
  X&=&\{Q:\Gamma\to SL(n,\C)\text{ analytic }\}=\{Q:\calO_\epsilon^{+}\cap\calO_\epsilon^{-}\to SL(n,\C) \text{ holomorphic }\},\\
  X^+_c&=&\{E:\calO_\epsilon^{+}\to SL(n,\C), \text{ holomorphic for some } \epsilon>0,\\ 
  &&\hspace{3in} E(p)=1, \text{ for some } p\in\calO_\epsilon \},\\
  X^-&=&\{F:\calO_\epsilon^{-}\to SL(n,\C), \text{ holomorphic for some }  \epsilon>0 \}.
\end{eqnarray*}
The Riemann-Hilbert problem is to factor $Q\in X$ into a product
$Q=E\cdot F,$ where $E\in X_c^+, F\in X^-$. This can be done on a big cell.

\begin{thm}
  $(X^+_c,X^-,X)$ is a local Manin triple of groups.
\end{thm}
The usual Riemann-Hilbert factorization scheme is given by the
following choices:
\begin{eqnarray}
  \Gamma&=&\{\lambda: |\lambda|=1\}\nonumber\\
  \calO^+_\epsilon&=&\{\lambda:|\lambda^{-1}|\leq 1+\epsilon\}\nonumber\\ 
  \calO^-_\epsilon&=&\{\lambda:|\lambda|\leq 1+\epsilon\}\label{3}
\end{eqnarray}
with $p=\infty\in\calO^+_\epsilon$ chosen as the normalizing point.

It is useful to keep in mind that the projection operators at the Lie
algebra level are given by Cauchy integral formulae. For $V\in\x=\{V:S^1\to
sl(n,\C) \text{ analytic }\}$
$$
\Pi^+V\in\x^+=\{W:\calO^+_\epsilon\to sl(n,\C) \text{ holomorphic }, W(\infty)=0\}.
$$
For $|\lambda|<1$,
\begin{equation}
  \label{4}
  \Pi^+V(\lambda)=W(\lambda)=\frac{1}{2\pi i}\oint_{|\lambda|=1}\frac{V(\xi)}{\lambda-\xi}\,d\xi. 
\end{equation}
Note that $\Pi^+V$ extends to a neighborhood of $\Gamma$ when $V$ is analytic
on $|\lambda|=1$.

We now procede to describe the more complicated Riemann-Hilbert
problem useful for harmonic maps. Now we have
\begin{eqnarray}
  \Gamma_\epsilon &=& \{ |\lambda|=\epsilon\}\cup\{|\lambda^{-1}|=\epsilon\}\nonumber\\
  \calO^+ &=& \C-\{0\}\nonumber\\
  \calO_\epsilon^- &=& \{|\lambda|<\epsilon\}\cup\{|\lambda^{-1}|\leq \epsilon\}.
\end{eqnarray}
The normalization point is chosen as $1\in\calO^+$. We also have a
reality condition.
\begin{definition}
  A map $Q:\calO\to SL(n,\C)$ satisfies the {\em harmonic map reality
    condition (HMRC)} if
  \begin{itemize}
    \item[(a)] $\lambda\to\bar{\lambda}^{-1}$ maps $\calO\to\calO$,
    \item[(b)] $Q(\lambda)=\left(Q(\bar{\lambda}^{-1})^*\right)^{-1}$.
  \end{itemize}
\end{definition}
This condition is compatible with the domains $\calO^+, \calO^-_\epsilon$ and
$\Gamma_\epsilon$, as well as the notion of holomorphicity.  Let
\begin{eqnarray}
  \label{6}
  Y^- &=& \{F:\calO^-_\epsilon\to SL(n,\C) \text{ holomorphic, } F \text{ satisfies the HMRC }\}\\
  Y^+ &=& \{E:\C-\{0\} \to SL(n,\C) \text{ holomorphic, } E \text{ satisfies the HMRC }\}\nonumber\\
  Y &=& \{Q:\calO^-_\epsilon -\left(\{0\}\cup\{\infty\}\right) \to SL(n,\C) \text{ holomorphic, }\nonumber 
    Q \text{ satisfies the HMRC }\}
\end{eqnarray}
As a special case of factorization, we have
\begin{thm}
  $(Y^+,Y^-,Y)$ is a local Manin triple.
\end{thm}

In the infinitesimal version, we have Lie algebras $(\y^+,\y^-,\y)$,
where $SL(n,\C)$ is replaced by $sl(n,\C)$. The reality condition
becomes
$$
W(\lambda)=-W(\bar{\lambda}^{-1})^*.
$$
The formula for the projection at the Lie algebra level is, for
$W\in\y$, where 
$$
\y=\{W:\calO^-_\epsilon - \left(\{0\} \cup\{\infty\}\right) \to sl(n,\C),
W(\lambda)=-W(\bar{\lambda}^{-1})^*\}
$$
For $\delta<\lambda<\delta^{-1}, \delta<\epsilon$
\begin{equation}
  \Pi^+W(\lambda)=\frac{1}{2\pi i}\int_{|\gamma|=\delta\cup|\gamma|=\delta^{-1}}\frac{W(\lambda)(\lambda-1)}{(\lambda-\gamma)(\gamma-1)}\,d\gamma 
\end{equation}
Since $W(\xi)$ is holomorphic in $\calO^-_\epsilon - \left(\{0\} \cup\{\infty\}\right)$,
the domain of analyticity of $\Pi^+W$ extends to $\C-\{\infty\}$.


\section{Derived Group Actions and (Half) Virasoro Actions}
\label{sec:4}

We now turn to a useful observation on passing automorphisms of $X$,
which restrict to automorphisms of $X^-$, to diffeomorphisms of $X^+$.
Assume $(X^+,X^-,X)$ is a local triple of groups.
\begin{definition}
  We call $G\subset \Hom(X,X)$ a negative automorphism group if
  $G\big|_{X^-}\subset \Hom(X^-,X^-)$. For $f\in G$, define $f^\# \in \Diff X^+$
  by
  \begin{equation}
    \label{8}
    f^\#(s_+)=P^+f(s_+).
  \end{equation}
\end{definition}

\begin{thm}
  $\#: G\to\Diff X^+$ is a (local) homomorphism.
\end{thm}
\begin{proof}
  Note that by local we mean $f^\#g^\# =(f\cdot g)^\#$, where $f^\#$ and $g^\#$
  are defined.We expect $f^\#$ and $g^\#$ to be well-defined close to
  $1\in G$.
  \begin{eqnarray*}
    (g\circ f)(s_+) &=& (g\circ f)^\#(s_+)r_1\\
    f(s_+) &=& f^\#(s_+)r_2\\
    g\left(f^\#(s_+)\right) &=& g^\#\left(f^\#(s_+)\right)r_3,
  \end{eqnarray*}
  where $r_1,r_2,r_3\in X^-$. Since $g$ is a group homorphism,
  \begin{eqnarray*}
    g\left(f(s_+)\right) &=& g\left(f^\#(s_+)r_2\right)\\
    &=& g\left(f^\#(s_+)\right)g(r_2)\\
    &=& g^\#\left(f^\#(s_+)\right)r_3g(r_2).\\
  \end{eqnarray*}
  By uniqueness, $r_1=r_3g(r_2)$ and
  $g^\#\left(f^\#(s_+)\right)=(g\circ f)^\#(s_+)$.
\end{proof}
\begin{cor}
  Let $\g$ and $(\x^+,\x^-,\x)$ be the Lie algebras of $G$ and
  $(X^+,X^-,X)$. Then the infinitesimal generator of the group action
  for $\sigma\in\g$ is given by
  \begin{equation}
    \label{9}
    \sigma^\#(s_+)=s_+\Pi^+\left(s_+^{-1}\sigma(s_+)\right),
  \end{equation}
  where 
  $$
  \#:X^+\times\g\to T(X^+).
  $$
  Moreover, $\#$ maps Lie brackets of elements in $\g$ to the Lie
  brackets of the vector-fields of the image.
\end{cor}
\begin{proof}
  This is the infinitesimal version of the group action. Hence the Lie
  bracket formula is a consequence of the composition law. To see the
  corectness of the formula for the infinitesimal generator, choose
  $\sigma\in\g$. Hence, for small $t$,
  $$
  e^{\sigma t}s_+=s_+(t)s_-(t)
  $$
  where $\left(s_+(t),s_-(t)\right)=\mu\inv\left(e^{\sigma t}s_+\right)$. But
  $$
  \sigma(s_+)=\sigma^\#(s_+)s_-(0) + s_+(0)\frac{d}{dt}\bigg|_{t=0} s_-(t).
  $$
  Since $s_-(0)=1$ and $\displaystyle\frac{d}{dt}\bigg|_{t=0}
  s_-(t)=V_-\in\x^-$.
  Notice also that $s_+(0)=s_+$, so we have
  $$
  s_+\inv\sigma(s_+)= s_+\inv\sigma^\#(s_+) + V_-,
  $$
  which implies
  $$
  \Pi^+\left(s_+\inv\sigma(s_+)\right)=s_+\inv\sigma^\#(s_+).
  $$
  Finally,
  $$
  \sigma^\#(s_+)=s_+\Pi^+\left(s_+\inv\sigma(s_+)\right).
  $$
  Note that one may start with this formula and prove directly that
  $\#$ is consistent with the Lie brackets if one wishes.
\end{proof}

The familiar example of a derived group action is the dressing action
of $X^-$ on $X^+$. Here the group action is
\begin{eqnarray*}
  \Ad X &\subseteq&\Hom(X,X)\\
  \Ad X^- &\subseteq&\Hom(X^+,X^+).
\end{eqnarray*}
The formula for the dressing action of $s_-\in X^-$ on $s_+\in X^-$ is
\begin{eqnarray*}
  s_-^\#s_+ &=& P^+(s_-s_+)=P^+(s_-s_+s_-\inv)\\
  &=& P_+\left(\Ad s_-(s_+)\right).
\end{eqnarray*}
The action can be thought of as derived either from left
multiplication or an $\Ad$ action. For more on dressing actions, see
either the book of Guest~\cite{G1} or the lecture notes of
Terng~\cite{T}.

We turn now to Virasoro actions. We persist in describing the actions
at the group level for conceptual simplicity, although we are
ultimately interested in infinitesimal formulae. The Virasoro algebra
is described in terms of generators
$$
V=\spann\{\dots,L_{-j},\dots, L_{-1},L_0,L_1,\dots\} 
$$
with the bracket operation $[L_j,L_k]=(k-j)L_{j+k}$. These are
suggestively written as $L_j=\lambda^{j+1}\displaystyle\dd{\lambda}$.  The span
can be over the complex numbers, yielding $\V_\C$, or the real
numbers, giving $\V_\R$. The algebras we will discover are
half-Virasoro algebras:
\begin{eqnarray*}
  \V_\C^+ &=& \spann_\C \{L_{-1},L_0,\dots,L_j\}\\
  \V_{\C,0}^+ &=& \spann_\C \{L_0,\dots,L_j,\dots\}\\
  \V_\C^- &=& \spann_\C \{\dots,L_{-j},\dots,L_0,L_1\}\\
  \V_{\C,\infty}^- &=& \spann_\C \{\dots,L_{-j},\dots,L_0\}.
\end{eqnarray*}
We will meet the corresponding real Virasoro algebras $\V_\R^\pm$, etc.
in section~\ref{sec:5}.

The Virasoro algebra is not the Lie algebra for a Lie group. However,
it can be considered as infinitesimal generators of holomorphic
mappings, which have composition properties modulo difficulties in
keeping track of domains and images. Since we are interested in
holomorphic mappings which are near the identity mapping, it is
possible to keep track of this, although we will be sloppy about it.
Choose a large set $\calO=\{\lambda:|\lambda|\leq N\}$ and let
$$
G^+_{hol}=\{f:f(\lambda)=\lambda+v(\lambda), v \text{ small and holomorphic in } \calO\}.
$$
If $f,g\in G^+_{hol}$, $f\circ g$ is defined and holomorphic, with
possibly a smaller domain. Let
$$
G^+_{hol,0}=\{f\in G^+_{hol}, f(0)=0\}.
$$
Then $\V_\C^+$ is the infinitesimal algebra for $G^+_{hol}$, and
$\V_{\C,0}^+$ is the infinitesimal algebra for $G^+_{hol,0}$.

Recall the usual triple $(X^+_c,X^-,X)$ for the Riemann-Hilbert
splitting, with the domains given in~\eqref{3}.
\begin{lemma}
  $G^+_{hol}$ acts by composition as a negative family of
  automorphisms on $(X^+_c,X^-,X)$.
\end{lemma}
\begin{proof}
  This is straightforward, except there is a real problem, which we do
  not try to solve, of keeping track of domains. Note that we have the
  notation
  \begin{equation}
    L_jF=\DD{t}\bigg|_{t=0} F\left(\lambda+t\lambda^{j+1}\right)=\lambda^{j+1}\dd{\lambda}F
  \end{equation}
  is the infinitesimal action of generators on $X^-$ or $X$.
\end{proof}

\begin{cor}
  \label{4.5}
  $G^+_{hol}$ has a derived action on $X^+_c$, where for $f\in G^+_{hol}, E\in X^+_c$
  $$
  f^\#E=P^+(E\circ f).
  $$
  For $V=v\displaystyle\dd{\lambda}\in\g^+_{hol}$,
  $$
  V^\#E=E\Pi^+\left(E\inv v\dd{\lambda}E\right).
  $$
\end{cor}
\begin{proof}
  This follows from~\eqref{2.1},~\eqref{2}, and~\eqref{9}.  Note that
  we have the explicit formula from~\eqref{4}.
  \begin{equation}
    \label{11}
    V^\#E(\lambda)=\frac{1}{2\pi i}E\oint_{|\gamma|=1}\frac{E\inv(\gamma)v(\gamma)\dd{\gamma}E(\gamma)}{(\lambda-\gamma)}\,d\gamma 
  \end{equation}
\end{proof}

We are now ready to construct the more complicated version of the
Virasoro action used for harmonic maps. First we construct
$G^+_{hol,0}$ and $\V_{\C,0}^+$ on $(Y_c^+,Y^-,Y)$ defined
in~\eqref{6}.  The domains $\calO_\epsilon^-$ consist of two pieces, a small
neighborhood about $0$, $\{|\lambda|\leq\epsilon\}$ and a small neighborhood about $\infty$,
$\{|\lambda|\inv\leq\epsilon\}$. $G^+_{hol,0}$ acts on $\{|\lambda|\leq\epsilon\}$. We induce the action
on $\{|\lambda|\inv\leq\epsilon\}$ by
$f(\lambda)=\left(\overline{f\left(\bar{\lambda}\inv\right)}\right)\inv$.
Now the action of $f$ is compatible with the HMRC~\eqref{2.2}.

\begin{thm}
  $G^+_{hol,0}$ acts by composition as a negative family of
  automorphisms on $(Y_c^+,Y^-,Y)$.
\end{thm}

\begin{cor}
  $G^+_{hol,0}$ has a derived action on $Y_c^+$, where for $f\in G^+_{hol,0}$
  $$
  f^\#E=EP^+\left(E\inv(E\circ f)\right)
  $$
  and for $V\in\V_{hol,0}^+, V=v\displaystyle\dd{\lambda}$,
  $$
  V^\#E=E\Pi^+\left(E\inv v\dd{\lambda}E\right).
  $$
\end{cor}
\begin{proof}
  The proof is the same as~\eqref{2.1} and~\eqref{2},
  despite the difficulties of keeping track of domains.
  However, we haven't really finished, as we are interested in
  explicit formulas. $V=v\displaystyle\dd{\lambda}$ is the correct
  expression at $\lambda=0$. However, at $\lambda=\infty$, the transformation
  $v(\lambda)=-\overline{v\left(\bar{\lambda}\inv\right)}$ gives
  $\overline{v\left(\bar{\lambda}\inv\right)}\lambda^2 \displaystyle\dd{\lambda}$.
\end{proof}

\begin{prop}
  If $V=v(\lambda)\displaystyle\dd{\lambda}$, then
  \begin{eqnarray*}
    V^\#E(\lambda) &=& \frac{1}{2\pi i} E\left[\oint_{|\gamma|=\epsilon}
      \frac{E\inv(\gamma)v(\gamma)\dd{\gamma}E(\gamma)(\lambda-1)}{(\lambda-\gamma)(\gamma-1)}\,d\gamma\right.\\
      &+& \left. \oint_{|\gamma|=\epsilon\inv}
        \frac{E\inv(\gamma)\overline{v(\bar{\gamma}\inv)}\gamma^2\dd{\gamma}E(\gamma)(\lambda-1)}{(\lambda-\gamma)(\gamma-1)}\,d\gamma\right]      
  \end{eqnarray*}
\end{prop}
Notice that the change of variable gives a correspondence between
$\V_{\C,0}^+$ and $\V_{\C,\infty}^-$:
$$
\sum_{j=0}^\infty c_j\lambda^{j+1}\dd{\lambda} \to -\sum_{j=0}^\infty \bar{c}_j\lambda^{-j+1}\dd{\lambda}.
$$
We are, in fact, choosing the graph of this representation in
$\V_{\C,0}^+\times\V_{\C,\infty}^-$, where the first Virasoro factor acts at
$\lambda=0$ and the second at $\lambda=\infty$. The full algebra
$\V_{\C,0}^+\times\V_{\C,\infty}^-$ would act on harmonic maps $s:\Omega\to SL(n,\C)$,
i.e. maps without the reality condition $s^*=s\inv$.

We give the formula for the generators. Take note that constants
multiply the second formula by the complex conjugate, so the formula
is a bit misleading. For $j\geq 0$, we have
\begin{eqnarray}
  \label{12}
  L_j(E)(\lambda) &=& \frac{1}{2\pi i} \left[\oint_{|\gamma|=\epsilon}
    \frac{E\inv(\gamma)\gamma^{j+1}\dd{\gamma}E(\gamma)(\lambda-1)}{(\lambda-\gamma)(\gamma-1)}\,d\gamma\right.\nonumber\\
  &+& \left. \oint_{|\gamma|=\epsilon\inv}
    x\frac{E\inv(\gamma)\gamma^{-j+1}\dd{\gamma}E(\gamma)(\lambda-1)}{(\lambda-\gamma)(\gamma-1)}\,d\gamma\right]      
\end{eqnarray}
Note that the singularities in the contour integral are at
$(0,\infty,1,\lambda)$. Hence there are many deformations of the contour possible
if we wish to compute $L_j(E)(\lambda)$ for $|\lambda|=1$ only.


\section{Virasoro Actions on Harmonic Maps}
\label{sec:5}

We associate to every harmonic map $s:\Omega\to SU(n)$ the extended harmonic
map $E:\C-\{0\}\times\Omega\to SL(n,\C)$ as in Theorem~\ref{2.3}. The Virasoro
action acts on the extended harmonic map via its dependence on the
variable $\lambda\in\C-\{0\}$. The special variable $z=x+iy$ is carried along as
an auxiliary variable.

\begin{thm}
  \label{5.1}
  Let $f\in G^+_{hol,0}$ be a holomorphic map near the identity, which
  we extend to $\calO_\epsilon^-$ by
  $f(\lambda)=\overline{f\left(\bar{\lambda}\in\right)\inv}$. Define
  \begin{equation}
    \label{13}
    \widehat{E}=f^*E=P^+(E\circ f).
  \end{equation}
  Let $\widehat{\Omega}=\{z\in\Omega : f^\#E_{\bullet}(z) \text{ is defined }\}$. Then
  $\widehat{E}$ is an extended  harmonic map on $\widehat{\Omega}$.
\end{thm}
\begin{proof}
  We need to show that $\widehat{E}$ satisfies the conditions of
  Theorem~\ref{2.4}. We have used Corollary~\ref{4.5} to define
  $\widehat{E}$. 
  
  Certainly $\widehat{E_1}(z)=I$ and $\widehat{E_\lambda}(p)=I$ by
  construction. The HMRC ensures that
  $\widehat{E_\lambda}(z)=\left(\widehat{E}_{\bar{\lambda}\inv}(z)\right)^{*-1}$.
  It is sufficient to show that $\widehat{E}_\lambda\inv\dd{z}E_\lambda$ has a
  simple pole at $0$ (and no poles at $\infty$). Then
  $\widehat{E}_\lambda\inv\dd{z}E_\lambda=\alpha+\lambda\inv \beta$, but $\beta=-\alpha$ since
  $\widehat{E}_1\inv\dd{z}E_1=0$. The HMRC gives a relationship
  between $\widehat{E}_\lambda\inv\dd{z}E_\lambda$ and
  $\widehat{E}_\lambda\inv\dd{\bar{z}}E_\lambda$ which finishes the proof.
  
  First, by construction, $\widehat{E}_\lambda(z)\in Y_c^+$ is holomorphic in
  $\lambda\in\C-\{0\}$. Hence $\dd{z}\widehat{E}_\lambda$ and $\widehat{E}_\lambda\inv$ and
  their products are holomorphic for $\lambda\in\C-\{0\}$. We need only worry
  about the singularity at $0$ and $\infty$. To handle this case, we look
  at the factors
  $$
  E_{f(\lambda)}(z)=\widehat{E}_\lambda(z)R_\lambda(z),
  $$
  where $R_\lambda(z)\in Y^-$ will be smooth in $z$. Hence we can write
  $\widehat{E}_\lambda=E_{f(\lambda)}R_\lambda\inv$ and 
  $$
  \widehat{E}_\lambda\inv\dd{z}\widehat{E}_\lambda=R_\lambda\dd{z}R_\lambda\inv +
  R_\lambda E_{f(\lambda)}\inv\dd{z}E_{f(\lambda)}R_\lambda\inv.
  $$
  The terms $R_\lambda, R_\lambda\inv, \dd{z}R_\lambda\inv$ are all holomorphic at
  $0$ and $\infty$. Hence the singularities at $0$ and $\infty$ come from the
  term $E_{f(\lambda)}\inv\dd{z}E_{f(\lambda)}=\left(1-f(\lambda)\inv\right)\alpha$. This
  expression is holomorphic at $\infty$. At $\lambda=0$, since $f(0)=0$ and $f$
  is ``close to the identity'', we have $E_{f(\lambda)}\inv\dd{z}E_{f(\lambda)}$
  has a simple pole at $\lambda=0$. The result follows.
\end{proof}

\begin{cor}
  \label{5.2}
  Let $V=v\dd{\lambda}\in\V_{\C,0}^+$. Then
  \begin{eqnarray}
    \label{14}
    V^\#E_\lambda(z) &=& \frac{1}{2\pi i}E_\lambda(z)\left[\oint_{|\gamma|=\epsilon}
      \frac{E_\gamma(z)\inv v(\gamma)\dd{\gamma}E_\gamma(z)(\lambda-1)}{(\lambda-\gamma)(\gamma-1)}\,d\gamma\right.\nonumber\\
    &+& \left.\oint_{|\gamma|=\epsilon\inv}
      \frac{E_\gamma(z)\inv \overline{v(\bar{\gamma})}\dd{\gamma}E_\gamma(z)(\lambda-1)}{(\lambda-\gamma)(\gamma-1)}\,d\gamma\right]
  \end{eqnarray}
  is tangent to the space of extended harmonic maps.
\end{cor}
\begin{proof}
  While for a holomorphic map $f\in G_{hol,0}^+$, the factorization
  defining the new extended harmonic map cannot always be done, if
  $f(\lambda)=\lambda+tv(\lambda)$, the factorization can be done if $t$ is sufficiently
  small. Since $f_t^\#E$ is an extended harmonic map for small $t$,
  $$
  V^\#E=\DD{t}\bigg|_{t=0}f_t^\#E
  $$
  os tangent to the space of extended harmonic maps.
\end{proof}

\begin{cor}
  \label{5.3}
  For $V\in\V_{hol,0}^+$ the map $V\to V^\#$ given in~\eqref{14} is a
  representation of $\V_{hol,0}^+$ on vector fields tangent to
  extended harmonic maps.
\end{cor}
\begin{proof}
  This is an application of Corollary~\ref{4.5} and Corollary~\ref{5.2}.
\end{proof}

\begin{cor}
  \label{5.4}
  Let
  \begin{eqnarray*}
    L_j(E_\lambda)(z) &=& \frac{1}{2\pi i} E_\lambda(z)\left[\oint_{|\gamma|=\epsilon}
      \frac{E_\gamma\inv(z)\gamma^{j+1}\dd{\gamma}E_\gamma(z)(\lambda-1)}{(\lambda-\gamma)(\gamma-1)}\,d\gamma\right.\\
    &+& \left. \oint_{|\gamma|=\epsilon\inv}
      \frac{E_\gamma\inv(z)\gamma^{-j+1}\dd{\gamma}E_\gamma(z)(\lambda-1)}{(\lambda-\gamma)(\gamma-1)}\,d\gamma\right]      
  \end{eqnarray*}
  generates a representation of $\V_{hol,0}^+$ on the space of vector
  fields tangent to extended harmonic maps.
\end{cor}
\begin{proof}
  This is a result of Corollary~\ref{5.3} in terms of specific
  generators. Note multiplication by a constant acts via
  multiplication by itself on the first factor, and to complex
  conjugate on the second.
\end{proof}
Corollary~\ref{5.4} is the main result of the paper. The reader is
invited to prove it directly and to decide whether the route we have
taken sheds light on the formula.


\section{The Results of Schwarz}

We now turn to the Virasoro action on harmonic maps from $\R^{1,1}$ to
$SU(n)$ treated by John Schwarz~\cite{Sch1}. The results are stated
without proof, as the construction is identical, except for the
different reality conditions.

\begin{thm}
  Let $E:\C^*\times\R^{1,1}\to SL(n,\C)$ be holomorphic in $\lambda\in\C^*$ and
  smooth for $(\xi,\eta)\in\R^{1,1}$. Assume $E_\lambda(0)=I$ and
  \begin{itemize}
  \item[(a)] $E_\lambda=\left(E_{\bar{\lambda}}\inv\right)^*$.
  \item[(b)] $E_1=I$.
  \item[(c)] $E_\lambda\inv\dd{\xi}E_\lambda$ and $E_\lambda\inv\dd{\eta}E_\lambda$ have simple
    pols at $0$ and $\infty$, respectively.
  \end{itemize}
  Then $s=E_{-1}:\R^{1,1}\to SU(n)$ is harmonic. Moreover, any harmonic
  map $s:\R^{1,1}\to SU(n)$ with $s(0)=I$ has an unique extended
  harmonic map associated with it.
\end{thm}

The different reality condition $f(\lambda)=\overline{f(\bar{\lambda})}$ results
in the decoupling of the Virasoro action at $0$ and $\infty$, as well as a
new restriction that the Virasoro actions are real.

\begin{thm}
  Let $W=w\dd{\lambda}\in\V_{\R,0}^+$ and $V=v\dd{\lambda}\in\V_{\R,\infty}^-$ be
  the elements of the product of two half-Virasoro algebras. Then
  \begin{eqnarray*}
    \delta_{v,w}(E_\lambda) &=& \frac{1}{2\pi i} E_\lambda\left[\oint_{|\gamma|=\epsilon}
      \frac{E_\gamma\inv w(\gamma)\dd{\gamma}E_\gamma(\lambda-1)}{(\lambda-\gamma)(\gamma-1)}\,d\gamma\right.\\
    &+& \left. \oint_{|\gamma|=\epsilon\inv}
      \frac{E_\gamma\inv v(\gamma)\dd{\gamma}E_\gamma(z)(\lambda-1)}{(\lambda-\gamma)(\gamma-1)}\,d\gamma\right]      
  \end{eqnarray*}
  is a representation of $\V_{\R,0}^+\times \V_{\R,\infty}^-$ on the vector
  fields tangent to the space of extended harmonic maps.
\end{thm}

Now Schwarz's formulae are quite different from these formulae.
However, a simple transformation $t=\frac{\lambda-1}{\lambda+1}$ and
$\tau=\frac{\gamma-1}{\gamma+1}$ will transform our integrals into his integrals.
However, because he is working with a different complex parameter, the
natural choice of Virasoro generators are $L_j=t^{j+1}\dd{t}$, which
are in principle allowable, since the expression in $\lambda$ is holomorphic
at $t=\pm 1$, or $\lambda=0,\infty$.

\begin{prop}
  $\V_\R \subset \V_\R^+ \times \V_\R^-$.
\end{prop}
\begin{proof}
  This embedding is achieved by expressing the generators $t^j\dd{t}$
  in terms of the coordinates adapted to $\pm 1$, or
  $t=\frac{\lambda-1}{\lambda+1}$. The algebra $\V_\R^+ \times \V_\R^-$ is actually
  larger, as the linear fractional transformations $L_j, j=-1,0,1$
  correspond to 3 generators in $\V_\R$ and 6 generators in $\V_\R^+ \times
  \V_\R^-$.
\end{proof}

Unfortunately, Schwarz fails to obtain a representation of the full
Virasoro algebra for the same reasons that we fail. We generate
$\V_{\R,0}^+ \times \V_{\R,\infty}^-\subset \V_\R^+ \times \V_\R^-$. A careful check of the
conditions of Schwarz shows that he also imposes the constraint that
the vector fields fix $(0,\infty)\sim (1,-1)$ and hence miss a full
realization. If we transform the description to one in terms of
scattering data, the harmonic maps correspond to first flows. The
missing $L_{-1}\in\V_\R^+$ and $L_1\in\V_\R^-$ will be written in terms of
second flows.


\end{document}